\documentclass[11pt]{amsart}

\newtheorem{theorem}{Theorem}[section]

\theoremstyle{definition}

\newtheorem{example}[theorem]{Example}

\theoremstyle{remark}

\newtheorem{claim}{Claim}[section]

\theoremstyle{theoremA}
\newtheorem*{theoremA}{Theorem A}
\theoremstyle{theoremB}
\newtheorem*{theoremB}{Theorem B}  

\numberwithin{equation}{section}

\newcommand{\real}{\mathbb{R}}

\newcommand{\Sp}{\mathbb{S}}

\newcommand{\la}{\lambda}
\newcommand{\La}{\Lambda}

\newcommand{\om}{\omega}
\newcommand{\Om}{\Omega}
\newcommand{\na}{\nabla}

\newcommand{\ka}{\kappa}
\newcommand{\De}{\Delta}

\newcommand{\vol}{\text{\rm vol}}

\newcommand{\lan}{\left\langle}
\newcommand{\ran}{\right\rangle}
\newcommand{\p}{\partial}

\newcommand{\dv}{\mathrm{div}}

\begin{document}

\title{The non-parabolicity of infinite volume ends}
\author{M. P. Cavalcante}
\address{Instituto de Matem\'atica, Universidade Fe\-deral de Alagoas, Macei\'o, 
AL, CEP 57072-970, Brazil}\email{marcos.petrucio@pq.cnpq.br}
\thanks{The first and the third authors were partially supported by CNPq under the
grants 483268/2010-0}
\author{H. Mirandola}
\address{Instituto de Matem\'atica, Universidade Fe\-deral do Rio de Janeiro, 
Rio de Janeiro, RJ, CEP 21945-970, Brasil} \email{mirandola@im.ufrj.br}
\author{F. Vit\'orio}
\address{Instituto de Matem\'atica, Universidade Fe\-deral de Alagoas, Macei\'o, 
AL, CEP 57072-970, Brazil}\email{feliciano.vitorio@pq.cnpq.br}

\date{January 1, 2001 and, in revised form, June 22, 2001.}


\begin{abstract}
 Let $M^m$, with $m\geq 3$, be an $m$-dimensional complete 
noncompact  manifold isometrically immersed in a Hadamard manifold $\bar M$. 
Assume that the mean curvature vector has finite $L^p$-norm, for some $2\leq p\leq m$. 
We prove that each end of $M$ must either have finite volume or be non-parabolic.  
\end{abstract}

\maketitle


\section{Introduction}

Let $(M^m, \langle,\rangle )$ be a complete noncompact Riemannian manifold without boundary. 
{We recall that $M$ is \emph{parabolic} if it does not admit a non-constant positive superharmonic 
 function. Otherwise, it is said to be \emph{non-parabolic}. There exist equivalent definitions for parabolic manifolds (see for instance Theorem 5.1 of \cite{g})}.
Let $E\subset M$ be an \emph{end} of $M$, that is an unbounded connected component of
$M-\Om$, for some compact subset $\Om\subset M$.
The property of parabolicity can be localized on each end of $M$.
Namely, we say that an end $E$ is \emph{parabolic} (see Definition 2.4 of \cite{l4}) 
if it does not admit a harmonic function $f:E\to \real$ satisfying:
\begin{enumerate}
\item $f|_{\p E}=1$;
\item $\liminf_{{y\to \infty}\atop{y\in E}} f(y)<1$.
\end{enumerate}
Otherwise, we say that $E$ is a \emph{non-parabolic} end of $M$. 
{It is well known that $M$ is non-parabolic if and only if it admits a non-parabolic end}. Furthermore, ends with finite volume are parabolic (see for instance  Section 14.4 of 
\cite{g}). In this direction we recall the following result due to Li and Wang:
\begin{theoremA}[Corollary 4 of \cite{lw} and Corollary 2.9 of \cite{l4}] Let $E$ be an end 
of a complete manifold. Suppose that, for some constants $\nu \geq 1$ and $C>0$, $E$ 
satisfies a Sobolev-type inequality of the form
\begin{equation}\label{sob-nu}
\Big(\int_E |u|^{2\nu}\Big)^{\frac{1}{\nu}} \leq C\int_E |\na u|^2,
\end{equation} 
for all compactly supported Sobolev function $u \in  W_c^{1,2}(E)$. 
Then $E$ must either have finite volume or be non-parabolic. Moreover, in the case $\nu>1$,  
$E$ must be non-parabolic.
\end{theoremA}

{Note that if a complete manifold $M$ that satisfies a Sobolev inequality as in Theorem A with $\nu=1$ (that is just the Dirichlet Poincar\'e inequality) then the first eigenvalue $\la_1(M)$ of the Laplace-Beltrami operator is positive, hence $M$ must be non-parabolic (see Proposition 10.1 of \cite{g}).
Example \ref{finite_volume_end} below exhibits a complete manifold that contains a finite volume end and that also satisfies a Sobolev inequality as in Theorem A with $\nu=1$.}

Cao, Shen and Zhu \cite{csz} showed that if $M^m$, with $m\geq 3$, is a complete manifold 
then each end of $M$ is non-parabolic provided that $M$ can be realized as a minimal 
submanifold in a Euclidean space $\real^{n}$. 
The same conclusion also was obtained by Fu and Xu \cite{fx} provided that there exists an 
isometric immersion of $M$ in a Hadamard manifold $\bar M$ with finite total mean 
curvature, that is, the mean curvature vector field $H$ of the immersion satisfies 
$\|H\|_{L^m(M)}<\infty$. In the both cases, they observed that $M$ admits a 
Sobolev-type inequality as in Theorem A with $\nu>1$.  

Our main result states the following:
\begin{theorem}\label{nonparabolic} Let $x:M^m\to \bar M$, with $m\geq 3$, be an isometric immersion of a complete non-compact manifold $M$ in a Hadamard manifold $\bar M$. Let $E$ be an end of $M$ such that the mean curvature vector satisfies $\|H\|_{L^p(E)}<\infty$, for some $2\leq p\leq m$. Then $E$ must either have finite volume or be non-parabolic.  
\end{theorem}
Example \ref{exemplofinitevolume} below exhibits an example of a complete non-compact hypersurface $M^m$ in $\real^{m+1}$, with $m\geq 3$, of finite volume and  mean curvature vector with  finite $L^p$-norm, for all $2\le p < m-1$. This example shows that Theorem \ref{nonparabolic} is not a consequence of Theorem A (except when $p=m$). Note also that the catenoids in $\real^3$ are parabolic minimal surfaces whose ends have infinite area, which shows that the hypothesis $m\geq 3$ is essential. 

In the present paper we also give a unified proof of the following fact:
\begin{theoremB} Let $x: M \to \bar M$ be an isometric immersion of a complete non-compact manifold $M$ in a manifold $\bar M$ with bounded geometry (i.e., $\bar M$ has sectional curvature bounded from above and injectivity radius bounded from below by a positive constant). Let $E$ be an end of $M$ and assume that the mean curvature vector of $x$ satisfies $\|H\|_{L^p(E)}<\infty$, for some $m\leq p \leq \infty$. Then $E$ must have infinite volume.
\end{theoremB}
The fact above was proved by Frensel \cite{fr} and by do Carmo, Wang and Xia \cite{cwx} for the case that the mean curvature vector field is bounded in norm (the case $p=\infty$),  by Fu and Xu \cite{fx} for the case that the total mean curvature is finite (the case $p=m$) and  by Cheung and Leung \cite{cl} for the case that the mean curvature vector has finite $L^p$-norm for some $p>m$. Since the cylinders of the form $M^m=\Sp^{m-1}\times \real$, where $\Sp^{m-1}$ is the unit Euclidean $(m-1)$-dimensional sphere, are examples of complete parabolic hypersurfaces in $\real^{m+1}$ we conclude that boundedness of the mean curvature vector does not imply that $M$ admits a Sobolev-type inequality. Furthermore, for all $m\geq 3$, we exhibit an example of a parabolic complete noncompact hypersurface $M^m$ in $\real^{m+1}$ such that the mean curvature vector has finite $L^p$-norm, for all $p>2(m-1)$. These examples show that Theorem B is not a consequence of Theorem A. 

Two questions arise in this paper: is there an example of a complete noncompact submanifold $M^m$, with $m\geq 3$, in a Euclidean space satisfying one of the conditions below?
\begin{enumerate}
\item $M$ has finite volume and $\|H\|_{L^p(M)}<\infty$, for some $m-1\leq p<m$;
\item $M$ is parabolic and $\|H\|_{L^p(M)}<\infty$, for some $m<p\leq 2(m-1)$.
\end{enumerate}

\section{Proof of Theorem \ref{nonparabolic}}
Choose $r_0>0$ so that the geodesic ball $B_{r_0}\subset M$ of radius $r_0$ and center at some point $\xi_0\in M$ satisfies $\p E\subset B_{r_0}$. For each $r>r_0$, consider $E_r=E\cap B_r$ and let 
$f_r:\overline E_r\to \real$ be a solution of 
the Dirichlet Problem:
\begin{equation*}
\left\{
\begin{array}{ll}
\De_M f_r=0& \mbox{ in } E_r,\\
f_r=1& \mbox{ in } \p E,\\
f_r=0& \mbox{ on } E\cap \p B_r.
\end{array}\right.
\end{equation*}
It follows from the maximum principle that $0< f_r\leq f_s<1$ in $E_r$, for all $s\geq r$.
Hence, by standard gradient estimates it follows that $\{f_r\}$ is an 
equicontinous family which converges uniformly on compact subsets,
when $r$ goes to infinity,  to a function $f:E\to \real$ satisfying
\begin{equation*}\label{fs}
\left\{
\begin{array}{ll}
\De_M f=0& \mbox{ in } E,\\
0\leq f\leq 1& \mbox{ in } E,\\
f=1&\mbox{ on } \p E.  
\end{array}\right.
\end{equation*}
If $f\not\equiv 1$  then it follows from the maximum principle that 
$\liminf_{x\to E(\infty)}f(x)<1$, which shows that $E$ is nonparabolic. 
Furthermore, it is well known that an end of finite volume is parabolic 
(see section 14.4 of \cite{g}). 
Hence, to prove Theorem \ref{nonparabolic}, it is sufficient to show the following:

\begin{claim}\label{claimpar} Either $f\not\equiv 1$ or $\vol(E)$ is finite.
\end{claim}

Suppose, by contradiction, that $f\equiv 1$ and $\vol(E)$ is infinite. This implies that, given any $L>1$, there exists
$r_1>r_0$ such that $\vol(E_{r_1}-E_{r_0})>2L$.
Since  $f_r\to 1$  uniformly on compact subsets, there exists
$r_2>r_1$ such that $f_r^{\frac{2m}{m-2}}>\frac{1}{2}$ everywhere in 
$E_{r_1}$, for all $r>r_2$. Thus, defining  $h(r)=\int_{E_r-E_{r_0}} f_r^{\frac{2m}{m-2}}$, with $r>r_0$, we obtain
\begin{equation}\label{h(r)}
h(r)\geq \int_{E_{r_1}-E_{r_0}} f_r^{\frac{2m}{m-2}}>L,
\end{equation}
for all $r>r_2$. In particular, we have that $\lim_{r\to\infty}h(r)=\infty$.

Now, for each $r>r_0$,  let   
$\varphi=\varphi_r \in C^\infty_0(E)$ be a cut-off function satisfying: 
\begin{enumerate}
\item $0\leq \varphi \leq 1$ everywhere in $E$;
\item $\varphi \equiv 1$ in $ E_r-E_{r_0}$.
\end{enumerate}
By Hoffmann-Spruck Inequality \cite{hs} we have
\begin{eqnarray*}\label{hsvarphi}
S^{-1}\left(\int_{E_r} (\varphi f_r)^{\frac{2m}{m-2}}\right)^{\frac{m-2}{m}}
\leq \int_{E_r} |\na (\varphi f_r)|^2 + \int_{E_r} (\varphi f_r)^2|H|^2,
 \end{eqnarray*}
where $S$ is a positive constant.

Note that
$$|\na (\varphi f_r)|^2=f_r^2|\na \varphi|^2+\varphi^2|\na f_r|^2 + 
\frac{1}{2}\lan\na\varphi^2,\na f_r^2\ran$$ and 
$$\varphi^2|\na f_r|^2 =\dv_M((f_r\varphi^2)\na f_r)-\frac{1}{2}\lan\na\varphi^2,\na f_r^2\ran,$$
since $f_r$ is harmonic. Using that $f_r\varphi$ vanishes on $\partial E_r$ we obtain
\begin{eqnarray*}\label{hsvarphi}
S^{-1}\left(\int_{E_r} (\varphi f_r)^{\frac{2m}{m-2}}\right)^{\frac{m-2}{m}}
&\leq& \int_{E_r} f_r^2|\nabla \varphi|^2+ \int_{E_r} \dv_M((f_r\varphi^2)\na f_r)
 \\&&+ \int_{E_r} (\varphi f_r)^2|H|^2\\
&=& \int_{E_r} f_r^2|\nabla \varphi|^2 + \int_{E_r} (\varphi f_r)^2|H|^2.
 \end{eqnarray*}
Thus, since $0\leq\varphi\leq 1$ in $E$ and $\varphi\equiv 1$ in $E_r-E_{r_0}$, we obtain
\begin{equation}\label{sob-c}
S^{-1}h(r)^{\frac{m-2}{m}}\leq S^{-1}\left(\int_{E_r} (\varphi f_r)^{\frac{2m}{m-2}}\right)^{\frac{m-2}{m}}
\leq \int_{E_{r_0}} f_r^2|\na \varphi|^2 + \int_{E_r} f_r^2|H|^2.
\end{equation}

First, assume that $\|H\|_{L^2(E)}$ is finite. Then, since $0\leq f_r\leq 1$, we have
\begin{equation*}\label{sob-fRp=2}
S^{-1} h(r)^\frac{m-2}{m}
\leq \int_{E_{r_0}}|\na \varphi|^2 + \int_{M}|H|^2.
\end{equation*}
Thus, $\lim_{r\to\infty}h(r)<\infty$, which is a contradiction.
Now, assume that $\|H\|_{L^p(E)}$ is finite, for some  $2<p\leq m$. Note that $\frac{m}{m-2}\leq \frac{p}{p-2}$. Since $0\leq f_r\leq 1$ and $h(r)>1$, for all $r>r_2$, we have: 
\begin{enumerate}
\item\label{ia} $f_r^{\frac{2p}{p-2}}\leq f_r^{\frac{2m}{m-2}}$;
\item\label{ib} $h(r)^\frac{p-2}{p}\leq h(r)^\frac{m-2}{m}$, for all $r>r_2$. 
\end{enumerate}
Thus, using H\"older Inequality, we have
\begin{eqnarray}\label{hder}
\int_{E_r-E_{r_0}} f_r^2 |H|^2 &\leq& \|H\|^2_{L^p(E_r-E_{r_0})}\left(\int_{E_r-E_{r_0}} f_r^{\frac{2p}{p-2}}\right)^{\frac{p-2}{p}}\\&\leq& \|H\|^2_{L^p(E-E_{r_0})}h(r)^{\frac{m-2}{m}},\nonumber 
\end{eqnarray}
for all $r> r_2$.

Choose $r_0>0$ large so that $\|H\|^2_{L^p(E-E_{r_0})}<\frac{1}{2S}$. Using  (\ref{sob-c}) and (\ref{hder}) we obtain the following:
$$ S^{-1} h(r)^{\frac{m-2}{m}}
\leq \int_{E_{r_0}} |\na \varphi|^2+\int_{E_{r_0}}|H|^2 + \ \frac{S^{-1}}{2}h(r)^{\frac{m-2}{m}}.$$
This shows that $\lim_{r\to\infty}h(r)<\infty$, which is a contradiction.
Therefore, Claim \ref{claimpar} and Theorem \ref{nonparabolic} are proved.


\section{Proof of Theorem B}  
Since $\bar M$ has bounded geometry, the sectional curvature $\bar K$ and the injectivity radius $i(\bar M)$ of $\bar M$ satisfy: \begin{equation}
\bar K<b^2 \mbox{ and } i(\bar M)>i_0,
\end{equation} 
for some positive constants $b$ and $i_0$. 
Let $E$ be an end of $M$ and assume that $\|H\|_{L^p(E)}$ is finite, for some $m\leq p\leq \infty$. Fix $\xi_0\in M$ and consider the geodesic balls $B_R=B_R(\xi_0)$, for all $R>0$.

Assume, by contradiction, that $\vol(E)$ is finite. The Hoffman-Spruck's isoperimetric inequality (see Theorem 2.2 of \cite{hs}) states that there exist a positive constant $\La$ and $S$, depending only $b$, $i_0$ and $m=\dim M$, such that if $N$ is a compact domain of $M$ with Lipschitz boundary $\p N$ then it holds:
\begin{equation}\label{isop_ineq}
S^{-1}\vol(N)^{\frac{m-1}{m}} \le \vol(\p N) + \int_{N} |H|,
\end{equation}
provided that the volume $\vol(N)<\La$.
Take $R_0>0$ sufficiently large so that $\p E\subset B_{R_0}$ and  $\vol(E-B_{R_0})<\min\{\La, 1\}$. Since $\|H\|_{L^p(E)}$ is finite, for some $m\le p\le +\infty$, we can take $R_0$ sufficiently large to satisfy further:
\begin{equation}\label{finite_lp_norm}
\begin{array}{ll}
\|H\|_{L^p(E-B_{R_0})}<\frac{1}{2S};& \mbox{ if } m\le p<+\infty;\\ 
\|H\|_{L^\infty(E)}\vol(E-B_{R_0})^{\frac{1}{m}}<\frac{1}{2S};& \mbox{ if } p=+\infty.
\end{array}
\end{equation}

Fix $R_1>R_0$, sufficiently large,  so that the distance $d_M(\p E,x)>R_0$, for all $x$ in $E-B_{R_1}$. Given any $q\in E-B_{2R_1}$ and $0<R<R_1$ we have that $B_R(q)\subset E-B_{R_1}$. Thus, by H\"older inequality, we obtain
\begin{equation}\label{argumento}
\begin{array}{ll}
\int_{B_R(q)}|H| \le \|H\|_{L^p(E-B_{R_1})}\vol(B_R(q))^{\frac{p-1}{p}}; \mbox{ if } m\le p<\infty;\\
\int_{B_R(q)}|H| \le \|H\|_{L^\infty(E)}\vol(E-B_{R_0})^{\frac{1}{m}}\vol(B_R(q))^{\frac{m-1}{m}}; \mbox{ if } p=\infty.
\end{array}
\end{equation}
Moreover, if $m\le p<\infty$ then $\frac{p-1}{p}\ge \frac{m-1}{m}$, which implies 
\begin{equation}\label{menor1}
\vol(B_R(q))^{\frac{p-1}{p}}\le \vol(B_R(q))^{\frac{m-1}{m}},
\end{equation}
since $\vol(B_R(q))\le \vol(E-B_{R_0})<1$. 

Since the distance function of $M$ from $\xi_0$ is a Lipschitz function, by using (\ref{isop_ineq}), (\ref{finite_lp_norm}), (\ref{argumento}) and (\ref{menor1}), we obtain
\begin{equation}\label{meancurvatureestimates}
\frac{1}{2S}\vol(B_R(q))^{\frac{m-1}{m}} \le \vol(\p B_R(q)).
\end{equation}
By the coarea formula, we have that $\vol(\p B_R(q))=\frac{d}{dR}\vol(B_R(q))$. Thus using (\ref{meancurvatureestimates}) we obtain
$\frac{d}{dR}(\vol(B_R(q)))^{\frac{1}{m}}\geq \frac{1}{2Sm}$.
This implies that
\begin{equation}\label{estimate} \vol(B_R(q))\geq \frac{1}{2Sm}R^m,
\end{equation}
for all $q\in E-B_{R_1}$ and $0<R<R_1$.

Since $M$ is complete and $E\subset M$ is connected and unbounded, there exists a sequence $p_2, p_3, \ldots$ in $E$ such that 
\begin{equation*}
p_k\in E\cap \left(B_{2kR_1}-B_{(2k-1)R_1}\right).
\end{equation*}
Note that $B_{{R_0}}(p_k)\subset E-B_{R_1}$ and $B_{{R_0}}(p_k)\cap B_{{R_0}}(p_{k'})=\emptyset$, for all $k\neq k'$. Since 
\begin{equation*}
\vol(E)\geq \vol(E-B_{R_1})\geq \sum_{k=2}^\infty \vol(B_{{R_0}}(p_k)),
\end{equation*}
 it follows from (\ref{estimate}) that $\vol(E)$ is infinite, which is a contradiction.  Theorem B is proved.


\section{Examples}\label{sec_examples}
\begin{example}\label{finite_volume_end} Consider the warped product manifold $M^m=\real\times_{e^{t}}P$, where $P$ is any complete $(m-1)$- dimensional manifold with finite volume. The metric of $M$ is complete and the end $E=(-\infty,0)\times P\subset M$ has finite volume given by 
\begin{equation*}
\vol(E)=\int_{-\infty}^0\int_P e^{m-1}dtdP=\frac{\vol(P)}{m-1} . 
\end{equation*}

Fix $k\in \real$ and let $h_\ka:M\to \real$ be the function defined by $h_\ka(t,x)=\ka t$. The gradient vector field of $h_\ka$ satisfies  
\begin{equation}\label{grad}
\na h_\ka=\ka \frac{\p}{\p t},
\end{equation}
where $\frac{\p}{\p t}(t,x)=\frac{d}{ds}\bigm|_{s=t}(s,x)\in T_{(t,x)}M$.  It is 
simple to show that $\na_Z \frac{\p}{\p t}=Z-\lan Z,\frac{\p}{\p t}\ran \frac{\p}{\p t}$. This implies that the Laplacian of $h_\ka$ satisfies 
\begin{equation}\label{laplacian}
\De h_\ka = \ka\,\dv\big(\frac{\p}{\p t}\big)=\ka\,(m-1).
\end{equation}
Fix $\eta\in C_0^{\infty}(M)$. Using (\ref{grad}) and (\ref{laplacian}) we obtain
\begin{eqnarray*}
\ka(m-1)\int_M \eta^2&=& \int_M \eta^2 \De h_\ka = \int_M (\dv (\eta^2 \na h_\ka) - 2\eta \lan \na \eta,\na h_\ka\ran)\\ &=& -2 \int_M \lan \na \eta, \eta \na h_\ka\ran \geq  - \int_M |\na\eta|^2 - \eta^2|\na h_\ka|\\&=& -\int_M |\na\eta|^2 - \ka^2 |\eta|^2.
\end{eqnarray*}
Thus, it holds that 
$\int_M |\na \eta|^2 + \ka(\ka + (m-1)) \eta^2 \geq 0$, for all $k\in \real$. In particular, if we take $\ka=-\frac{m-1}{2}$ we obtain
\begin{equation*}
\int_M |\na\eta|^2 -\frac{(m-1)^2}{4}\eta^2 \geq 0. 
\end{equation*}
Hence $M$ satisfies a Sobolev inequality as in Theorem A with $\nu=1$.
\end{example}

\begin{example}\label{nonparabolic_end} Let $f:(-\infty,\infty)\to (0,\infty)$ be 
a positive smooth function satisfying that $f(t)=f(-t)$ and $f(t)=t^{\frac{1}{m-1}}$, 
for all $t\geq 1$. 
Consider the immersion $x:\Sp^{m-1}\times \real\to \real^{m+1}=\real^m\times \real$ 
given by $x(v,t)= (f(t)v, t)$.
Consider $M$ the product $\Sp^{m-1}\times \real$ endowed with the metric induced by $x$. 
The metric of $M$ is given by
\begin{equation}\label{metric_example_nonparabolic}
\lan\,,\ran_{(v,t)}=(1+f'(t)^2) dt^2 + f(t)^2 \lan\,,\ran_v, 
\end{equation}
where $\lan\,,\ran_v$ denotes the metric of $\Sp^{m-1}$. Note that $M$ is a complete manifold with two ends. 

We claim that $M$ is parabolic. To do this,  it is sufficient to prove that the following ends of $M$: 
\begin{equation*}
E_+=(1,\infty)\times \Sp^{m-1} \mbox{ and } E_-=(-\infty,-1)\times \Sp^{m-1}.
\end{equation*}
are parabolic (see Proposition 14.1 of \cite{g}).
In fact, we define:
\begin{equation*}
V_+(s)=\vol_M\left(\{q\in E_+ \bigm| d(q,\p E_+)\leq s\}\right)
\end{equation*}
and 
\begin{equation*}
V_-(s)=\vol_M\left(\{q\in E_- \bigm| d(q,\p E_-)\leq s\}\right)
\end{equation*}
Using (\ref{metric_example_nonparabolic}) and that 
$f(t)=t^{\frac{1}{m-1}}$, for all $|t|\geq 1$, we obtain that 
$$V_+(s)=V_-(s)\leq { D s^2},$$ for some constant $D>0$ and for all $s\geq 1$. 
In particular, $$\int^\infty \frac{s}{V_+(s)}ds=\int^\infty \frac{s}{V_+(s)}ds=\infty.$$ 
This implies that $M$ is parabolic (see section 14.4 of \cite{g}). 

We claim that the mean curvature vector $H$ of the isometric immersion $x$ has 
finite $L^p$-norm, for all $p>m$. In fact, a simple computation shows that
\begin{equation}
m H({x(v,t)})= \frac{(m-1)}{f(t)\sqrt{1+f'(t)^2}} - \frac{f''(t)}{(1+f'(t)^2)^{\frac{3}{2}}}
\end{equation}
Using that $f(t)=t^{\frac{1}{m-1}}$, for all $|t|\geq 1$, we obtain that 
$|H({x(v,t)})|\leq C t^{-\frac{1}{m-1}}$, for some $C>0$ and for all 
$x(v,t)\in E_+\cup E_-$. Thus, we obtain
\begin{equation*}
\int_M |H|^p dM \leq D \int^\infty t^{1-\frac{p}{m-1}}dt,
\end{equation*}
for some $D>0$. This implies that $\|H\|_{L^p(M)}$ is finite when $p>2(m-1)$.
\end{example}

\begin{example}\label{exemplofinitevolume} Let $x:\Sp^{m-1}\times \real\to \real^{m+1}$ be the 
immersion given by $x(v,t)= (e^{-t^2}v, t)$ and
consider $M=\Sp^{m-1}\times \real$ endowed with the metric induced by $x$. 
The metric of $M$ is complete and the volume of $M$ is given by
\begin{eqnarray}
\vol(M)&=&\om_{m-1}\int_{-\infty}^{+\infty} (1+4t^2e^{-2t^2})^{\frac{1}{2}}e^{-(m-1)t^2}dt,
\end{eqnarray}
where $\om_{m-1}$ is the volume of $\Sp^{m-1}$. This implies that $\vol(M)$ is finite, since 
the integral $\int_{-\infty}^{+\infty} e^{-(m-1)t^2}dt$ is finite and the function 
$t\in \real\mapsto 1+4t^2e^{-2t^2}$ is bounded. In particular, $M$ is parabolic 
since it has finite volume (see Theorem 7.3 of \cite{g}).

The mean curvature vector $H$ of the isometric immersion $x$ is given by 
\begin{equation}
H({x(v,t)})=h(t)=\frac{2e^{-t^2}(1-2t^2)}{m(4t^2e^{-2t^2}+1)^\frac{3}{2}} + 
\frac{(m-1)e^{t^2}}{m(4t^2e^{-2t^2}+1)^{\frac{1}{2}}}.
\end{equation}
Using that $\lim_{t\to\infty}e^{-t^2}(1-2t^2)=\lim_{t\to\infty}4t^2e^{-2t^2}=0$ we obtain that 
\begin{equation}
\lim_{t\to\pm\infty} h(t)e^{-t^2}= \frac{m-1}{m}.\end{equation}
Thus the integral
\begin{eqnarray}
\int_M |H|^p&=& \om_{n-1}\int_{-\infty}^{+\infty} \left(|h(t)|^p 
(1+4t^2e^{-2t^2})^{\frac{1}{2}}e^{-(m-1)t^2}\right)dt \nonumber\\
&=& \om_{n-1}\int_{-\infty}^{+\infty} \left((|h(t)|e^{-t^2})^p 
(1+4t^2e^{-2t^2})^{\frac{1}{2}}\right)e^{(p-m+1)t^2}dt
\end{eqnarray}
converges if $0\leq p<m-1$ and diverges if $p\geq m-1$.
\end{example}

\section*{Acknowledgement} The authors would like to thank Detang Zhou for helpful comments during the preparation of this article. The third named author thanks Universidade Federal do Rio de Janeiro for hospitality during preparation of this article.

\bibliographystyle{amsplain}


\end{document}